\documentstyle[12pt]{article}

\normalsize

\begin{document}
\baselineskip 17pt
\title{{\bf On the new identities of Dirichlet $L$-functions}
\footnotetext{This work is supported by Natural Science Basic Research Project of Shaanxi Province(2021JM-044).\\
1.marong0109@163.com\qquad 2.pi2cookie.zhang@gmail.com\qquad
 3.\,zzboyzyl@163.com }
\author{Rong Ma$^{1}$ \quad Jinglei Zhang$^{2}$ \\
\small{School of Mathematics and Statistics,
Northwestern Polytechnical University, Xi'an, Shaanxi} \\
\small{ 710072,
People's Republic of China}\\
Yulong Zhang$^{3}$\\ \small{School of Software Engineering, Xi'an Jiaotong University
Xi'an, Shaanxi}\\
\small{ 710049,
People's Republic of China}}
\date{}}
\maketitle{}
\vspace{0.0cm}
\begin{center}
\begin{minipage}{125mm}
\begin{center}
\large{{\bf Abstract}}
\end{center}
Let $q\ge3$ be an integer, $\chi$ be a Dirichlet character
modulo $q$, and $L(s,\chi)$ denote the Dirichlet $L$-functions
corresponding to $\chi$. In this paper, we show some special function series, and
give some new identities for the Dirichlet $L$-functions involving Gauss sums. Specially, we give specific identities for $L(2,\chi)$.
\end{minipage}
\end{center}

\noindent {\large\textbf{1.Introduction}}\\

Let $q\ge3$ be an integer and $\chi$ be a Dirichlet character modulo $q$. Dirichlet $L$-functions $L(s,\chi)$ is defined by
\[L(s,\chi)=\sum_{n=1}^{\infty}\frac{\chi(n)}{n^s},\]
where $s$ is a complex number with $Re(s)>1$. It can be
extended to all $s$ using analytic continuation. Many scholars have studied the mean value of Dirichlet $L$-functions and have got some identities or asymptotic formulae [1-10]. For example, Alkan [1] has got the identical equations of Dirichlet $L$-functions when $s=1,2$ with the odd and even Dirichlet character $\chi$ respectively,
\begin{equation}
	\sum_{\mbox{\tiny$\begin{array}{c}
\chi\bmod{q}\\ \chi(-1)=-1
\end{array}$}}|L(1,\chi)|^2=\frac{\pi^2\varphi(q)}{12}\prod_{p|q}\left(1-\frac{1}{p^2}\right)-\frac{\pi^2\varphi^2(q)}{4q^2};	
\end{equation}
\begin{equation}
\sum_{\mbox{\tiny$\begin{array}{c}
\chi\bmod{q}\\ \chi(-1)=1
\end{array}$}}|L(2,\chi)|^2=\frac{\pi^4\varphi(q)}{180}\prod_{p|q}\left(1-\frac{1}{p^4}\right)
+\frac{\pi^4\varphi(q)}{18q^2}\prod_{p|q}\left(1-\frac{1}{p^2}\right),
\end{equation}
where $\varphi(q)$ is the Euler's totient function.
In addition, in contrast with these results, he also has given approximate formulas which are relatively easy to be derived by Abel's summation and the Polya-Vinogradov inequality, that is
\begin{equation}
	\sum_{\mbox{\tiny$\begin{array}{c}
\chi\bmod{q}\\ \chi(-1)=-1
\end{array}$}}|L(1,\chi)|^2=\frac{\varphi(q)}{2}+O(\sqrt{q}\log q),
\end{equation}
where $\varphi(q)$ is defined the same as before and the signal $O$ only depends on constants.

He also has remarked that if $\chi$ and $r$ have the same parity, then it is always possible to determine the average value of $|L(r,\chi)|^2$ when $\chi$ ranges over all odd or even character modula $q$.

 To proceed, we shall introduce the Gauss sum $G(z,\chi)$ [3] corresponding to $\chi$ modulo $q$ and Bernoulli numbers [3] which are defined as
\[G(z,\chi):=\sum_{k=0}^{q-1}\chi(k)e^{\frac{2\pi ikz}{q}};\]
\[B_n=\sum_{k=0}^{n}{n \choose k}B_k,\]
where $z$ is an integer, ${n\choose k}=\frac{n!}{k!(n-k)!}$ is the binomial coefficient and $B_0=1$.

In the proof of formulae (1) and (2), Alkan [2] has mainly constructed an identity between the Gauss sum and Dirichlet $L$-functions which is satisfied for the special conditions, that is,
\begin{equation}
	\frac{(-1)^{s+1}qs!}{i^s2^{s-1}\pi^s}L(s,\chi)=\sum_{j=1}^{q}G(j,\chi)
\sum_{k=0}^{2[\frac{s}{2}]}{s \choose k}B_k\left(\frac{j}{q}\right)^{s-k},
\end{equation}
where $s\ge 1$ and $\chi$ have the same parity.

Of course, Zhang [14] also has obtained formulae (1) and (2) with the help of the identity between the Dirichlet $L$-functions and the generalized Dedekind sum
\begin{eqnarray}
S(h,n,q)&=&\frac{(n!)^2}{4^{n-1}\pi^{2n}q^{2n-1}}\sum_{d|q}\frac{d^{2n}}{\varphi(d)}
\sum_{\tiny\begin{array}{c}
\chi\bmod{q}\\ \chi(-1)=(1)^{n}
\end{array}}\chi(h)|L(n,\chi)|^2\nonumber\\
&&-\frac{(n!)^2}{4^{n}\pi^{2n}}\left(\sum_{r=1}^{\infty}\frac{1+(-1)^n}{r^n}\right)^2,
\end{eqnarray}
where the generalized Dedekind sum is defined by
\[S(h,n,q)=\sum_{a=1}^{q}\bar{B}_n(\frac{a}{q})\bar{B}_n(\frac{ah}{q}),\]
which $\bar{B}_n(x)$ is defined as following
\begin{equation}
	\bar{B}_n(x)=
	\left\{\begin{array}{ll}
			B_n(x-[x]),&\mbox{if $x$ is not an integer,}\\
			0,&\mbox{if $x$ is an integer.}
		\end{array}\right.
\end{equation}

Although Alkan has obtained formulae (1) and (2) with the help of the connection in [2] between the Dirichlet L-function and Gauss sum, he can only give the connection for $s$ and $\chi$ having the same parity. In fact, when $s$ and $\chi$ had the different parity, Alkan has got the right of formula (4) is zero. Therefore following his method, he could not get the identities when $s$ and $\chi$ have the different parity. Based on that, we would like to study the connection between Dirichlet $L$-functions and Gauss sums under the condition that $s$ and $\chi$ have the different parity. By the method, we can also give a new connection between Dirichlet $L$-functions and Gauss sums under the condition $s$ and $\chi$ having the same parity. I think it's very interesting and significant because we can know more about the connections between values of Dirichlet L-functions on positive integers and Gauss sum and do more research about Dirichlet $L$-functons.

In this paper, we will generate identities between Dirichlet L-functions and Gauss sums, and get some new identities of Dirichlet L-functions. That is, we will prove the following:

\noindent \textbf{Theorem 1}
	Let $q\ge3$ be an integer and $\chi$ be a Dirichlet character modulo $q$. If $s>1$ is an integer, $s$ and $\chi$ have the different parity, we have
\begin{equation}
L(s,\chi)=\frac{(-i)^{s(\bmod 2)}}{q}\sum_{j=1}^{q}G(j,\chi)\sum_{n=0}^{\infty}\frac{\sin^{(s+1)(\bmod 2)} (\frac{2\pi n j}{q})\cos^{s(\bmod 2)}(\frac{2\pi n j}{q})}{n^s},
\end{equation}
where $G(j,\chi)$ is the Gauss sum corresponding the Dirichlet character $\chi$ modulo $q$.

\noindent \textbf{Theorem 2}
Let $q\ge3$ be an integer and $\chi$ be a Dirichlet character modulo $q$. If $s$ is an integer, $s$ and $\chi$ have the same parity, we have
\begin{equation}
L(s,\chi)=-\frac{(-i)^s}{qs!}\sum_{j=1}^qG(j,\chi)\left[\frac{x^s}{2}-\frac{s\pi x^{s-1}}{2}+\sum_{k=1}^{[ s/2]}\zeta(2k)(x^s)^{(2k)}\right]_{x=\frac{2\pi j}{q}},
\end{equation}
where $G(j,\chi)$ is the Gauss sum corresponding the Dirichlet character $\chi$ modulo $q$, $\zeta(2k)$ is the Riemann zeta function at value $2k$, and $[x]$ denotes the biggest integer less than $x$.

Specially, if $s=2$ and $\chi$ is the even or odd Dirichlet character modula $q$ respectively, we have

\noindent \textbf{Corollary 1}
	Let $q\ge3$ be an integer and $\chi$ be an even and odd Dirichlet character modulo $q$ respectively, we have
\begin{eqnarray}
L(2,\chi)&=&\frac{\pi^2}{q^3}\sum_{j=1}^{q}j^2G(j,\chi)-\frac{\pi^2}{q^2}\sum_{j=1}^{q}jG(j,\chi)+\frac{\pi^2}{6q},\chi(-1)=1;\\
L(2,\chi)&=&\frac{1}{q}\sum_{j=1}^{q}G(j,\chi)\int_{0}^{\frac{2\pi j}{q}}\log\left(2\sin\frac{x}{2}\right)dx,\chi(-1)=-1,
\end{eqnarray}
where $G(j,\chi)$ is the Gauss sum corresponding the Dirichlet character $\chi$ modulo $q$.

\textbf{Note.} Theoretically speaking, we can give all the identities for $L(r,\chi)$ for every integer $r\ge1$ for the corresponding even and odd Dirichlet character $\chi$ respectively except $L(1,\chi)$ for the even character because of the condition of Theorem 2 . How to get the identity of $L(1,\chi)$ for an even character is still an open problem.

\noindent {\large\textbf{2. Some Lemmas}}

To prove the theorems, we give the following lemmas.

\noindent \textbf{Lemma 1}
	If a function sequence $\{G_n(x)\}_{n=0,1,2\cdots}$ satisfies the differential equation

\[\frac{d}{dx}G_{n+1}(x)=(n+1)G_{n}(x),\]
with $G_0(x)$ is a nonzero integer. Then the function sequence will be determined by a constant sequence $\{a_n\}$ in the form of
	\[G_n(x)=\sum_{k=0}^{n}{n\choose k}a_kx^{n-k},\]
where $a_n$ is the constant term of polynomial $G_n(x)$.

\noindent \textbf{Proof:}
Obviously $G_n(x)$ is a polynomial of degree $n$, assume that \[G_n(x)=\sum_{k=0}^n{n \choose k}a_kx^{n-k}.\]
	
Multiplying both sides by $n+1$ and integrating, then we have
	\begin{eqnarray}
			\int{(n+1)G_{n}(x)dx}&=&(n+1)\sum_{k=0}^n{n \choose k}a_k\frac{x^{n-k+1}}{n-k+1}\nonumber\\
			&=&\sum_{k=0}^n\frac{(n+1)!}{k!(n-k+1)!}a_kx^{n-k+1}+a_{n+1}\nonumber\\
			&=&\sum_{k=0}^{n+1}{n+1\choose k}a_kx^{n+1-k}\nonumber\\
			&=&G_{n+1}(x).
	\end{eqnarray}
Taking the derivative of formula (11), then we can get
\begin{eqnarray}
			(n+1)G_{n}(x)=\frac{d}{dx}G_{n+1}(x).\nonumber
	\end{eqnarray}
This has proved Lemma 1.

\noindent \textbf{Lemma 2}
	If $x\in (0,2\pi)$, we have
\begin{eqnarray}	
&&\sum_{n=1}^{\infty}\frac{\sin nx}{n}=\frac{\pi-x}{2};\\
&&\sum_{n=1}^{\infty}\frac{\cos nx}{n}=-\log2\left(\sin \frac{x}{2}\right).
\end{eqnarray}

\noindent \textbf{Proof:} As we know, for a complex variable $z$, when $|z|<1$, we have the Taylor series expansion of $\log(1-z)$ as
\[-\log(1-z)=\sum_{n=1}^{\infty}\frac{z^n}{n}.\]

We can also extend it to all $|z|\le 1$ except $z=1$ by using analytic continuation. Therefore, let $z=e^{ix}$, it is allowed for $x\in (0,2\pi)$
\begin{equation}
-\log (1-e^{ix})=\sum_{n=1}^{\infty}\frac{\cos nx}{n}+i\sum_{n=1}^{\infty}\frac{\sin nx}{n}\;.
\end{equation}
Substitute $1-e^{ix}=\sqrt{2(1-\cos x)}e^{\frac{x-\pi}{2}i}$ into equation (14), we get
\begin{eqnarray}
		-\log (1-e^{ix})&=&-\log\sqrt{2(1-\cos x)}-\frac{x-\pi}{2}i\nonumber\\
		&=&-\log\left(2\sin\frac{x}{2}\right)-\frac{x-\pi}{2}i\nonumber\\
		&=&\sum_{n=1}^{\infty}\frac{\cos nx}{n}+i\sum_{n=1}^{\infty}\frac{\sin nx}{n}\nonumber
\end{eqnarray}

Consider real part and imaginary part respectively, we immediately have
\begin{eqnarray}
	&&\sum_{n=1}^{\infty}\frac{\sin nx}{n}=\frac{\pi-x}{2};\nonumber\\
	&&\sum_{n=1}^{\infty}\frac{\cos nx}{n}=-\log2\left(\sin \frac{x}{2}\right).\nonumber
\end{eqnarray}

That proves Lemma 2.

\noindent {\large\textbf{2. Proofs of Theorems}}

With the help of these lemmas, it is feasible to prove theorems and the corollary. Firstly we will prove Theorem 1.

\textbf{\noindent Proof of Theorem 1.} In order to prove Theorem 1, we need to discuss two equations between Dirichlet $L$-functions and Gauss sums under the different condition $\chi(-1)=1$ and $\chi(-1)=-1$  respectively.

When $\chi(-1)=1$ and $s>1$, we have
\begin{eqnarray}%
	&&\sum_{j=1}^qG(j,\chi)\sum_{n=1}^\infty\frac{\cos \frac{2\pi nj}{q}}{n^s}\nonumber\\
&=&\sum_{j=1}^q\sum_{m=1}^q\chi(m)e^{\frac{2\pi mji}{q}}\sum_{n=1}^\infty\frac{e^{\frac{2\pi nji}{q}}+e^{\frac{-2\pi nji}{q}}}{2n^s}\nonumber\\
&=&\sum_{n=1}^\infty\frac{1}{2n}\sum_{m=1}^q\chi(m)\sum_{j=1}^q\left(e^{\frac{2\pi (m+n)ji}{q}}+e^{\frac{2\pi (m-n)ji}{q}}\right)\nonumber\\
&=&\sum_{n=1}^\infty\frac{q}{2n^s}[\chi(-n)+\chi(n)]\nonumber\\
&=&q\sum_{n=1}^\infty\frac{\chi(n)}{n^s}\nonumber\\
&=&qL(s,\chi),\nonumber
\end{eqnarray}
which is equivalent to
\begin{equation}
	L(s,\chi)=\frac{1}{q}\sum_{j=1}^qG(j,\chi)\sum_{n=1}^\infty\frac{\sin (\frac{2\pi n j}{q})}{n^s}.
\end{equation}

When $\chi(-1)=-1$, we have
\begin{eqnarray}
		&&i\sum_{j=1}^qG(j,\chi)\sum_{n=1}^\infty \frac{\sin (\frac{2\pi n j}{q})}{n^s}\nonumber\\
&=&\sum_{j=1}^q\sum_{m=1}^q\chi(m)e^{\frac{2\pi mji}{q}}\sum_{n=1}^\infty\frac{e^{\frac{2\pi nji}{q}}-e^{\frac{-2\pi nji}{q}}}{2n^s}\nonumber\\
&=&\sum_{n=1}^\infty\frac{1}{2n^s}\sum_{m=1}^q\chi(m)\sum_{j=1}^q\left(e^{\frac{2\pi (m+n)ji}{q}}-e^{\frac{2\pi (m-n)ji}{q}}\right)\nonumber\\
&=&\sum_{n=1}^\infty\frac{q}{2n^s}[\chi(-n)-\chi(n)]\nonumber\\
&=&-q\sum_{n=1}^\infty\frac{\chi(n)}{n^s}\nonumber\\
&=&-qL(s,\chi),\nonumber
\end{eqnarray}
which is equivalent to
\begin{equation}
	L(s,\chi)=-\frac{i}{q}\sum_{j=1}^qG(j,\chi)\sum_{n=1}^\infty\frac{\cos (\frac{2\pi n j}{q})}{n^s}.
\end{equation}

According to the equations (15) and (16), when $s$ and $\chi$ have the same parity, we have
\begin{equation}
	L(s,\chi)=\frac{(-i)^{s(\bmod\;2)}}{q}\sum_{j=1}^qG(j,\chi)\sum_{n=1}^{\infty}\frac{\sin^{s(\bmod\;2)} (\frac{2\pi n j}{q})\cos^{(s+1)(\bmod\;2)} (\frac{2\pi n j}{q})}{n^s}.
\end{equation}

When $s$ and $\chi$ have the different parity, $s>1$, we have
\[L(s,\chi)=\frac{(-i)^{s(\bmod\;2)}}{q}\sum_{j=1}^qG(j,\chi)\sum_{n=0}^{\infty}\frac{\sin^{(s+1)(\bmod\;2)} (\frac{2\pi j}{q})\cos^{s(\bmod\;2)} (\frac{2\pi j}{q})}{n^s}.\]

This proves Theorem 1. Next we will prove Theorem 2 from the formula (17).

\textbf{Proof of Theorem 2.} Consider the part of trigonometric function separately
\[\sum_{n=1}^\infty\frac{\sin^{s(\bmod\;2)} (\frac{2\pi n j}{q})\cos^{(s+1)(\bmod\;2)} (\frac{2\pi n j}{q})}{n^s},\]
which is equivalent to the following function sequence we denote as $\{F_s(x)\}$ at the point $x=\frac{2\pi j}{q}$
\begin{equation}
	\sum_{n=1}^{\infty}\frac{\sin nx}{n},\;\sum_{n=1}^{\infty}\frac{\cos nx}{n^2},\;\sum_{n=1}^{\infty}\frac{\sin nx}{n^3},\;\sum_{n=1}^{\infty}\frac{\cos nx}{n^4},\;\sum_{n=1}^{\infty}\frac{\sin nx}{n^5}\cdots
\end{equation}

Therefore we can rewrite the formula (17) as following
\begin{equation}
L(s,\chi)=\frac{(-i)^{s(\bmod\;2)}}{q}\sum_{j=1}^qG(j,\chi)F_s\left(\frac{2\pi j}{q}\right).
\end{equation}

For formula (18), we have recursion obviously,
\begin{eqnarray}
	\left\{\begin{array}{c}
		F_1(x)=\frac{\pi-x}{2}\\
		F_n(x)=(-1)^nF_{n+1}'(x),
	\end{array}\right.
\end{eqnarray}
where the leading term is obtained by Lemma 2. If we give another function sequence $\{G_n(x)\}$ which is defined by
\begin{eqnarray}
	\left\{\begin{array}{c}
		G_n(x)=-n!F_n(x),  \mbox{ when $n\equiv 0,1(\bmod 4)$};\\
		G_n(x)=n!F_n(x),  \mbox{ when $n\equiv 2,3 (\bmod 4)$}.
	\end{array}\right.
\end{eqnarray}
Then we have $\{G_n(x)\}$ is a function sequence satifying the requirement of Lemma 1. Therefore, a sequence $\{a_n\}$ exists to determine $\{G_n(x)\}$. What means, if we find the general term formula of $\{a_n\}$, we find that of $\{F_n(x)\}$.
When $0<x<2\pi$, through the formulae (18) and (20), we have
\begin{eqnarray}
		G_0(x)&=&1/2;\nonumber\\
		G_1(x)&=&\frac{x}{2}-\frac{\pi}{2};\nonumber\\
		G_2(x)&=&\frac{x^2}{2}-\pi x+\frac{\pi^2}{3};\nonumber\\
		G_3(x)&=&\frac{x^3}{2}-\frac{3\pi x^2}{2}+\pi^2x+0;\nonumber\\
		G_4(x)&=&\frac{x^4}{2}-2\pi x^3+\pi^2x^2+0x+\frac{4\pi^2}{15};\nonumber\\
		G_5(x)&=&\frac{x^5}{2}-\frac{5\pi x^4}{2}+\frac{10\pi^2x^3}{3}-\frac{4\pi^2x}{3}+0;\nonumber\\
		\cdots\cdots\nonumber
\end{eqnarray}
Then list the constant term of the polynomial $G_n(x)$, we get the sequence $\{a_n\}$ as following
\begin{eqnarray}
a_0&=&1/2;\nonumber\\
a_1&=&-\pi/2;\nonumber\\
a_{2k}&=&(2k)!\sum_{j=1}^\infty\frac{1}{n^{2k}};\nonumber\\
a_{2k+1}&=&0, \nonumber\\
&&(k=1,2,3\cdots)
\end{eqnarray}
where $\sum_{j=1}^\infty\frac{1}{n^{2k}}$ is the value of Riemann $\zeta$ function on the even integers as \[\zeta(2k)=\frac{(-1)^{k+1}B_{2k}(2\pi)^2k}{2(2k)!}\;.\]

Let $[t]$ denote the round down of $t$ and $(x^n)^{(t)}$ denote the $t$-order derivative of $x^n$, from Lemma 1 and formula (22), so we have
\begin{eqnarray}
		G_n(x)&=&\sum_{k=0}^{n}{n\choose k}a_kx^{n-k}\nonumber\\
		&=&\frac{x^n}{2}-\frac{n\pi x^{n-1}}{2}+\sum_{k=1}^{[ n/2]}(2k)!{n\choose 2k}\zeta(2k)x^{n-2k}\nonumber\\
		&=&\frac{x^n}{2}-\frac{n\pi x^{n-1}}{2}+\sum_{k=1}^{[ n/2]}\zeta(2k)(x^n)^{(2k)}.
\end{eqnarray}

Plug the formula (23) into the formula (21), considering the formula (19) and we get
\begin{eqnarray}
	&&L(s,\chi)=\nonumber\\
	&&\left\{\begin{array}{c}
		-\frac{1}{qs!}\sum_{j=1}^qG(j,\chi)\left[\frac{x^s}{2}-\frac{s\pi
x^{s-1}}{2}+\sum_{k=1}^{[s/2]}\zeta(2k)(x^s)^{(2k)}\right]_{x=\frac{2\pi j}{q}},\,
 \mbox{$s\equiv 0(\bmod 4)$};\\
		\frac{i}{qs!}\sum_{j=1}^qG(j,\chi)\left[\frac{x^s}{2}-\frac{s\pi x^{s-1}}{2}+\sum_{k=1}^{[s/2]}\zeta(2k)(x^s)^{(2k)}_{x=\frac{2\pi j}{q}}\right],\quad\quad\, \mbox{$s\equiv1\;(\bmod\;4)$};\\
		\frac{1}{qs!}\sum_{j=1}^qG(j,\chi)\left[\frac{x^s}{2}-\frac{s\pi x^{s-1}}{2}+\sum_{k=1}^{[s/2]}\zeta(2k)(x^s)^{(2k)}\right]_{x=\frac{2\pi j}{q}} ,\quad\, \mbox{$s\equiv2\;(\bmod\;4)$};\\
		-\frac{i}{qs!}\sum_{j=1}^qG(j,\chi)\left[\frac{x^s}{2}-\frac{s\pi x^{s-1}}{2}+\sum_{k=1}^{[s/2]}\zeta(2k)(x^s)^{(2k)}\right]_{x=\frac{2\pi j}{q}}, \, \mbox{$s\equiv3\;(\bmod\;4)$}.\nonumber
	\end{array}\right.
\end{eqnarray}
In conclusion
\[L(s,\chi)=-\frac{(-i)^s}{qs!}\sum_{j=1}^qG(j,\chi)\left[\frac{x^s}{2}-\frac{s\pi x^{s-1}}{2}+\sum_{k=1}^{[s/2]}\zeta(2k)(x^s)^{(2k)}\right]_{x=\frac{2\pi j}{q}}.\]

This proves Theorem 2.

\noindent {\textbf{Proof of Corollary 1}
We can easily get formula (9) from Theorem 2 for $s=2$. For formula (10), we remark that definitive integral appearing in the Corollary 1 is a improper integral. It converges because\[\int_0^\frac{\pi}{4}\log(2\sin\frac{x}{2})dx=-4\pi\log2.\]

Through the Lemma 2, when $x\in(0,2\pi)$
\begin{eqnarray*}
		\sum_{n=1}^{\infty}\frac{\sin nx}{n^2}=\int\sum_{n=1}^{\infty}\frac{\cos nx}{n}+C\\
		=\int_0^x-\log2(\sin \frac{t}{2})dt\\
\end{eqnarray*}
So when $s=2$ in Theorem 2, taking $x=\frac{2\pi j}{q}$, we have
\[L(2,\chi)=\frac{1}{q}\sum_{j=1}^q G(j,\chi)\int_0^{\frac{2\pi j}{q}}\log(2\sin\frac{x}{2})dx.\]

The proof of Corollary 1 has already been completed.

\end{document}